\documentclass[10pt]{article}

\usepackage{authblk}
\usepackage{amsmath}
\usepackage{amsfonts}
\usepackage{amssymb}
\usepackage{mathrsfs}
\usepackage{nicefrac}
\usepackage{mathtools}
\usepackage{subcaption}
\usepackage{booktabs}
\usepackage{graphicx}
\usepackage[top=2cm,
            bottom=2cm,
            right=1.5cm,
            left=1.5cm,
]{geometry}
\usepackage[linesnumbered]{algorithm2e}

\SetKwComment{Comment}{/* }{ */}
\RestyleAlgo{ruled}
\SetKwData{Input}{Input}
\DontPrintSemicolon

\providecommand{\keywords}[1]
{
  \small	
  \textbf{\textit{Keywords---}} #1
}

\begin{document}

\title{Hyperparameter optimization of orthogonal functions in the numerical solution of differential equations}

\author[1]{Alireza Afzal Aghaei\thanks{alirezaafzalaghaei@gmail.com}}

\author[1,2]{Kourosh Parand\thanks{k\_parand@sbu.ac.ir}}

\affil[1]{Department of Computer and Data Science, Faculty of Mathematical Science , Shahid Beheshti University, Tehran, Iran.}
\affil[2]{Department of Cognitive Modeling, Institute for Cognitive and Brain Sciences , Shahid Beheshti University, Tehran, Iran.}

\maketitle

\abstract{
This paper considers the hyperparameter optimization problem of mathematical techniques that arise in the numerical solution of differential and integral equations. The well-known approaches grid and random search, in a parallel algorithm manner, are developed to find the optimal set of hyperparameters. Employing rational Jacobi functions, we ran these algorithms on two nonlinear benchmark differential equations on the semi-infinite domain. The configurations contain different rational mappings along with their length scale parameter and the Jacobi functions parameters. These trials are configured on the collocation Least-Squares Support Vector Regression (CLS-SVR), a novel numerical simulation approach based on spectral methods. In addition, we have addressed the sensitivity of these hyperparameters on the numerical stability and convergence of the CLS-SVR model. The experiments show that this technique can effectively improve state-of-the-art results. 
 }


\section{Introduction}\label{sec1}
Estimating the unknown dynamics of a given physical system is an essential task in science and engineering. Mathematicians simulate these systems after expressing them in functional equations such as differential and integral equations. The failure of the analytical approaches in solving these problems, which mainly contain nonlinear terms, has led researchers to develop various numerical techniques. Finite Elements (FEM), Finite Volume (FVM), Meshless, and Spectral methods are some well-known efforts. These techniques consider a simple mathematical model defined by some parameters and hyperparameters. The former are internal variables learned during the training process, whereas the latter are external configurations that should be optimized to find the best estimator. These terms are initially back in the machine learning literature. For example, in a supervised machine learning task, the Support Vector Machine (SVM) algorithm considers a hyperplane for separating the given data into different classes.
The vector that defines the hyperplane is a parameter, whereas the kernel function is a hyperparameter. Likewise, a linear combination of unknown weights and basis functions is usually considered for approximating the solution to a differential equation. Here the unknown weights are parameters, while the basis functions are a hyperparameter. Choosing the best combination of hyperparameters is critical and may significantly affect the result. Scientists take advantage of prior knowledge to choose sub-optimal ones. As an example, it is common to use logarithmic-spaced values for the regularization parameter of a machine learning model. In a mathematics setting, the intrinsic nature of the problem is utilized to choose appropriate hyperparameters. For instance, the periodic problems are usually approximated by the Fourier basis functions, and the rational functions may approximate the problems that reach their steady state, and so on.

Regarding the development of artificial intelligence, hyperparameter optimization or simply tuning and finding the optimal hyperparameters, has evolved. Various approaches and techniques are proposed to handle this issue. Exhaustive search \cite{bergstra2012random}, probabilistic algorithms \cite{snoek2012practical}, gradient-based \cite{maclaurin2015gradient} and meta-heuristic algorithms \cite{itano2018extending} are some of these efforts. In a similar manner, mathematicians tried to develop some routines to find optimal hyperparameters that arises in the numerical simulation techniques. To name a few, Boyd  \cite{boyd1982optimization} proved some theorems for rational Chebyshev functions, Sanyasiraju et al. \cite{sanyasiraju2013optimization} used a local optimization algorithm for Radial Basis Functions (RBFs), Cavoretto et al. \cite{cavoretto2021search} combined leave-one-out cross-validation with univariate global optimization techniques for RBFs, Tanguy et al. \cite{tanguy2015parameter} presented a quadratic minimization procedure for optimal placement of poles in rational approximations by Müntz–Laguerre functions and Mi at al. \cite{9096566} developed two algorithms for the adaptive selection of continuous rational orthogonal basis functions. However, it is valuable to provide methods for solving problems in a broader range rather than just specific cases.

A wide range of physical problems involving initial and boundary value problems are defined on the infinite or semi-infinite domain. Approximating an accurate solution on the corresponding domain is essential. These problems are mostly solved with orthogonal polynomials or rational functions with completeness and orthogonality properties on the problem domain. The Hermite, Laguerre, and rational Jacobi functions are the most widely-used functions applied to these problems. Computational and mathematical properties of the rational functions encouraged scientists to choose them as the basis functions \cite{tajvidi2008modified,khader2020numerical,saadatmandi2018collocation, abd2021matrix,deniz2020rational, zhang2019revisiting, parand2011collocation, dehghan2015rational, parand2018generalized}. However, these functions suffer of hyperparameters such as rational mapping and length scale parameter. These can disturb the numerical solution in some cases \cite{boyd2001chebyshev}.

In this research, we develop some algorithms and investigate the applications of machine learning algorithms to optimize the hyperparameters that appear during the numerical simulation of differential equations arising on semi-infinite domains. In the continuation of the article, we will discuss the preliminaries (section 2), the proposed method (section 3), and the state-of-the-art numerical results (section 4). Finally, the concluding remarks will be discussed in the last section.

 \section{Preliminaries}\label{sec2}
In this section, we explain some prerequisites needed in the following sections. To do so, we first explain the Jacobi polynomials, then the CLS-SVR method will be recalled. The hyperparameter optimization techniques used in the rest of the work will be discussed.

\subsection{Orthogonal Polynomials}\label{subsec21}
Hermite, Laguerre, and Jacobi polynomials are the most well-known orthogonal polynomials defined on infinite, semi-infinite, and finite intervals. They can be used to approximate functions in corresponding domains. However, approximating rational functions may not be very accurate using polynomials. Therefore, researchers proposed various rational mappings to handle this issue. In general, Jacobi polynomials with hyperparameters $\alpha, \beta$ are defined on the interval $[-1, 1]$ by the recursive expression
\begin{equation}
\begin{aligned}
&J_{i}^{\alpha, \beta}(x)=-\frac{(\alpha+i-1)(\beta+i-1)(\alpha+\beta+2 i)}{i(\alpha+\beta+i)(\alpha+\beta+2 i-2)} J_{i-2}^{\alpha, \beta}(x) \\
&+\frac{(\alpha+\beta+2 i-1)\left\{\alpha^{2}-\beta^{2}+x(\alpha+\beta+2 i)(\alpha+\beta+2 i-2)\right\}}{2 i(\alpha+\beta+i)(\alpha+\beta+2 i-2)} \\
&\quad \times J_{i-1}^{\alpha, \beta}(x), \quad i=2,3, \ldots,
\end{aligned}
\label{eq:jacobi}
\end{equation}
where
\begin{equation*}
J_{0}^{\alpha, \beta}(x)=1, \quad J_{1}^{\alpha, \beta}(x)=\frac{\alpha+\beta+2}{2} x+\frac{\alpha-\beta}{2}
\end{equation*}
Their orthogonality is defined using the $L2$ inner product:
\begin{equation*}
    {\displaystyle \int _{-1}^{1}(1-x)^{\alpha }(1+x)^{\beta }J_{m}^{(\alpha ,\beta )}(x)J_{n}^{(\alpha ,\beta )}(x)\,dx={\frac {2^{\alpha +\beta +1}}{2n+\alpha +\beta +1}}{\frac {\Gamma (n+\alpha +1)\Gamma (n+\beta +1)}{\Gamma (n+\alpha +\beta +1)n!}}\delta _{nm}}.
\end{equation*}
The Gegenbauer, Chebyshev, and Legendre polynomials are some special cases of Jacobi polynomials. For Legendre polynomials, the equation \eqref{eq:jacobi} with $\alpha=\beta=0$ reduces to:
\begin{equation}
	\begin{aligned}
	\label{eq:leg}
	&P_0 (x)=1,\quad P_1 (x)=x, \\
	&(n+1) P_{n+1}(x)=(2n+1)xP_n(x)-nP_{n-1}(x),\quad n\ge 1,
	\end{aligned}
	\end{equation}
and for Chebyshev with $\alpha=\beta=-\nicefrac{1}{2}$ we have 
\begin{equation}
\begin{aligned}
&T_0 (x)=1,\quad T_1 (x)=x, \\
&T_{n+1}(x)=2x\,T_{n}(x)-T_{n-1}(x),\quad n\ge 1.\end{aligned}
\label{eq:cheb}
\end{equation}
Although these polynomials are defined for the bounded domain, researchers have used some nonlinear maps $\phi$ with the property $\phi:[-1,1]\to[0,\infty)$ to transform the orthogonality into the semi-infinite domain. To our best knowledge, these three mappings are the most used functions among researchers \cite{parand2018accurate}:
\begin{itemize}
    \item Algebraic mapping: \quad \qquad
    $\displaystyle\phi(x)=\nicefrac{(x-\theta)}{(x+\theta)}$.
    \item Exponential mapping: \qquad $\displaystyle\phi(x)=1-2{\textrm{e}^ {-\displaystyle\nicefrac{x}{\theta}}}$.
    \item Logarithmic mapping: \qquad $\displaystyle \phi(x)=2\tanh\left(\nicefrac{x}{\theta}\right)-1$.
\end{itemize}
The rational Jacobi functions is defined by a simple transformation $\tilde{J}(x) = J(\phi(x))$ where $\phi(x)$ is a nonlinear rational mapping. Therefore, orthogonality takes the form
\begin{equation}
    {\displaystyle \int _{0}^{\infty}\tilde{J}_{m}^{(\alpha ,\beta )}(x)\tilde{J}_{n}^{(\alpha ,\beta )}(x) w(x) \,dx={\frac {2^{\alpha +\beta +1}}{2n+\alpha +\beta +1}}{\frac {\Gamma (n+\alpha +1)\Gamma (n+\beta +1)}{\Gamma (n+\alpha +\beta +1)n!}}\delta _{nm}},
\end{equation}
where $w(x)$ is the corresponding weight function of the inner product.
\subsection{Collocation Least-Squares Support Vector Regression}\label{subsec22}
Collocation Least-Squares Support Vector Regression (CLS-SVR) is a novel formulation of support vector machines for solving functional equations \cite{khoee2022least, ahadian2022support, rad2023learning}. In this machine learning model, the unknown solution of a differential equation is approximated by a linear combination of unknown coefficients and some known basis functions. The basis functions, also known as feature maps, transform the input data into a nonlinear space in which we hope the approximation accuracy will increase. In the following, we recall the formulation of CLS-SVR which is based on the paper \cite{parand2021parallel}.

Suppose that we want to solve an arbitrary differential equation in the form of $\mathcal{N}(u) = f$. Approximating the solution by $m$ basis functions, we have:
\begin{equation*}
	\tilde{u}(x) = w^T\varphi(x) + b =  \sum_{i=1}^{m}w_{i}\varphi_{i}(x) + b.
\end{equation*}
The primal optimization form of CLS-SVR takes the form:
\begin{equation}
\begin{aligned}
	\min\limits_{w, e} \quad& \frac{1}{2} w^Tw + \frac{\gamma}{2}e^Te \\
	\text{ s.t.} \quad& \mathcal{N}(\tilde{u})(x_k) - f(x_k) = e_k, \quad k=1, \ldots , n,\\
	& \tilde{u}(c_j) = u_j, \quad j=1,\ldots,d.
    \label{eq:opt}
\end{aligned} 
\end{equation}
where $n$ is the number of training points, $e_k$ is the value of the residual function at $k$-th training point, and $\tilde{u}(c_j) = u_j$ is the set of initial and boundary conditions. The regularization hyperparameter $\gamma$ controls the fluctuations of the learned solution and reduces the overfitting on the training data. However, this hyperparameter may not be optimized in the case of problems with unique solutions.

The dual form of this optimization problem leads to a linear or nonlinear system of equations. This system can be obtained using the Lagrangian function for \eqref{eq:opt}:
\begin{align*}
    \mathscr{L}(w, e, \alpha) = \frac{1}{2} w^Tw + \frac{\gamma}{2}e^Te - \sum_{k=1}^n \alpha_k \left[\mathcal{N}(\tilde{u})(x_k) - f(x_k) - e_k \right] - \sum_{j=1}^d \beta_j \left[ \tilde{u}(c_i) - u_i \right].
	\end{align*}
The saddle point of this function satisfies the solution of the problem:
    \begin{equation*}
    \left\{\frac{\partial\mathscr{L}}{\partial w_i} = 0,
\frac{\partial\mathscr{L}}{\partial e_k} = 0,
\frac{\partial\mathscr{L}}{\partial\alpha_k} = 0, \frac{\partial\mathscr{L}}{\partial\beta_j} = 0\right\},
\end{equation*}
where $i=1,\cdots, m$, $j=1,\cdots,d$, and $k=1,\cdots, n$. After solving this system, we can use the obtained weights $w$ to approximate the equation or use the dual variables in a kernel method sense. 

\subsection{Hyperparameter optimization}\label{subsec23}
\begin{figure}[!htbp]
\centering
\includegraphics[width=.7\textwidth]{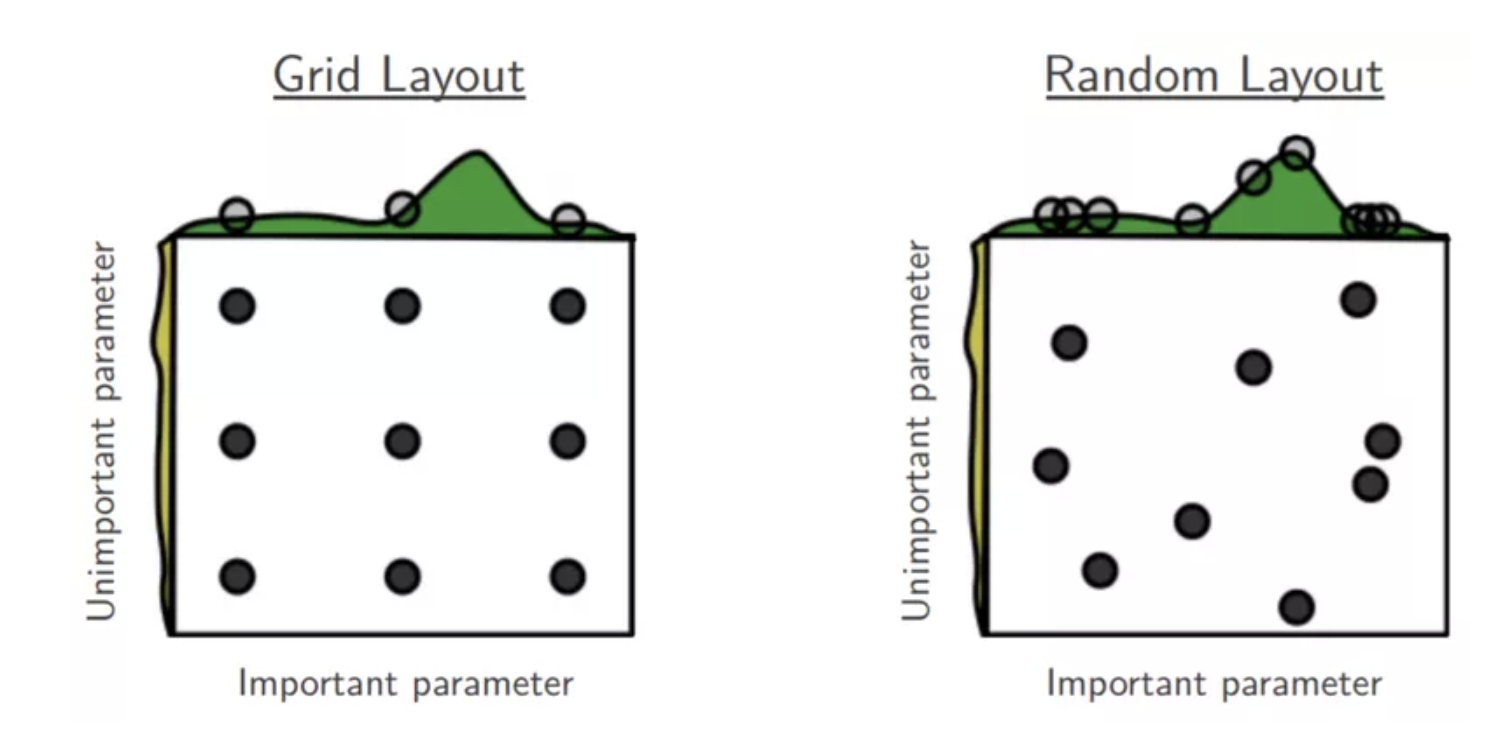}
\caption{. 1: A comparison between random and grid searches for a simple search space with only two different parameters. It can be seen that grid search may fail to explore the search space efficiently. The figure is adapted from Bergstra et al. \cite{bergstra2012random}. }\label{fig:gridvsrand}
\end{figure}

The problem of hyperparameter optimization is usually expressed as a non-convex minimization or maximization problem. Various algorithms proposed to solve this problem can be categorized into gradient-based and gradient-free sets. The former uses the gradient of a loss function to find an optimum value. At the same time, by generating some candidates, the latter tries to approximate the search space to find a local or possibly global optimum. Grid search, random search, Bayesian optimization, and meta-heuristic algorithms are some examples of gradient-free methods. Grid and random search have been widely used because of their simplicity and acceptable performance.
Moreover, they can be massively run in parallel. Briefly, the grid search seeks the Cartesian product of given parameter sets, whereas the random search samples a fixed number of parameters from user-defined distributions. For categorical cases, the hyperparameter is chosen uniformly. In both methods, the best parameter set would be measured on a test dataset usually generated by the K-Fold cross-validation algorithm. Figure \ref{fig:gridvsrand} compares these two algorithms. 

A key difference between grid and random search is their computational complexity. The grid search seeks all possible values in $O(n^k)$ time, having $k$ hyperparameters and $n$ different values for each of them, while for a random search, the user can define a budget based on the available resources so that the algorithm will run on the time complexity of $O(n)$. For a more precise comparison, we refer to the authors \cite{yang2020hyperparameter, yu2020hyper}.

\section{Method and Results}\label{sec3}

In this section, we explain the proposed algorithm and then provide examples to show the methods’ efficiency. Here, we focus on the hyperparameter optimization of orthogonal rational Jacobi functions, whereas the presented algorithm can be easily extended to other mathematical models. 

The optimal set of hyperparameters should be found on a test dataset to prevent the overfitting problem. In machine learning, this is easily done using cross-validation techniques. Nevertheless, in the mathematical equations, there is no data to split into train and test sets. However, there are three alternative options to handle this issue. The first option is to use a set of new nodes (different from training points) in the domain and calculate the prediction error on these nodes. The second is employing some physical properties of the model. The last one is to use some criteria found by mathematical techniques, including analytical or accurate numerical ones. The last option can be seen as a composition of previous ones. Most of the criteria measures used by the authors have important physical meanings which are good enough to be accuracy test criteria. Therefore, we use the absolute difference value between the exact or the state-of-the-art approximations and the predicted one given by CLS-SVR. The proposed grid and random search algorithms are presented in algorithms \ref{alg:proposedgrid} and \ref{alg:proposedrand}. 

We will obtain the optimal numerical results for some well-known benchmark nonlinear differential equations using the proposed method in the following sections. The configuration used to find an accurate numerical solution for these problems is reported in tables \ref{table:gridvalues} and \ref{table:randvalues}. It is seen that the grid search seeks $600$ different configurations to find the optimal set. To have a fair comparison, we set the maximum iterations for random search the same as the number of grid nodes. In addition, the roots of Legendre and Chebyshev polynomials are utilized as the training data.

\begin{center}
\begin{minipage}{0.45\textwidth}
\begin{algorithm}[H]
\KwData{The differential equation \textbf{As} $ODE$}
% \KwData{Accuracy criteria \textbf{As} $Criteria$}
$S_i \gets \text{Set of all desired values for } i\text{-th hyperparameter}$\;
$Criteria List \gets []$\;
$Parameter\_Set \gets CartesianProduct(\{S_i: \forall i\})$\;
\For{$param\_set \textbf{ in } Parameter\_Set \textit{ in parallel}$}{
$\tilde{u} \gets CLS-SVR(ODE)$\;
Compute Criteria for $\tilde{u}$\;
Push Criteria to $Criteria List$\;
}
\KwResult{$parameter\_set$ associated with  best criteria}
\caption{Grid search algorithm}\label{alg:proposedgrid}
\end{algorithm}
\end{minipage}
\hfill
\begin{minipage}{0.45\textwidth}
\begin{algorithm}[H]
\KwData{The differential equation \textbf{As} $ODE$}
% \KwData{Accuracy criteria \textbf{As} $Criteria$}
$S_i \gets \text{A suitable distribution for } i\text{-th hyperparameter}$\;
$Criteria List \gets []$\;
\;
\For{$iter \textbf{ from } 1 \textbf{ to } MAX\_ITER \textit{ in parallel}$}{
$\tilde{u} \gets CLS-SVR(ODE)$\;
Compute Criteria for $\tilde{u}$\;
Push Criteria to $Criteria List$\;
}
\KwResult{$parameter\_set$ associated with best criteria}
\caption{Random search algorithm}\label{alg:proposedrand}
\end{algorithm}
\end{minipage}
\end{center}

\begin{table}[!htbp]
\centering
\begin{minipage}{0.49\textwidth}
\centering
\begin{tabular}{@{}ll@{}}
\toprule
Parameter & Values \\ \midrule
Kernel & \{Legendre, Chebyshev\} \\
Mapping & \{Algebraic, Exponential, Logarithmic\} \\
$\theta$ & $\{0.1,0.2,\cdots, 9.9, 10\}$ \\ \bottomrule
\end{tabular}
\caption{The search space of the grid search} \label{table:gridvalues}
\end{minipage}
\begin{minipage}{0.49\textwidth}
\centering
\begin{tabular}{@{}ll@{}}
\toprule
Parameter & Values \\ \midrule
Kernel & \{Legendre, Chebyshev\} \\
Mapping & \{Algebraic, Exponential, Logarithmic\} \\
$\theta$ & UniformDistribution(0, 10) \\ \bottomrule
\end{tabular}
\caption{The search space of the random search} \label{table:randvalues}
\end{minipage}
\end{table}

\subsection{Volterra's population model}\label{subsec31}
Volterra's population model is a nonlinear integro-differential equation that describes population growth of a species within a closed system \cite{parand2021parallel}. It is common to express this problem into an equivalent nonlinear differential equation \cite{parand2011collocation}:
\begin{equation}
\begin{aligned}
&\kappa{u''}(x)={u'}(x) - {u'}(x)^2 - {u}(x) {u'}(x),\\
&{u}(0)= 0,
{u'}(0)=0.1.
\end{aligned}
\label{eq:volterrapopulation}
\end{equation}
Here $\kappa$ is a non-dimensional parameter. The criteria for the prediction correctness of this problem is the maximum value of the approximated function. TeBeest \cite{tebeest1997classroom} showed that the maximum peak is:
\begin{align}
u_{max} = 1 + \kappa\ln({\frac{\kappa}{1+\kappa+{u(0)}}}).
\label{eq:volterrapopulationexact}
\end{align}
Considering \eqref{eq:volterrapopulationexact} as the exact value for \eqref{eq:volterrapopulation}, the absolute error for this equations is defined as
\begin{align}
\text{Absolute Error}\vcentcolon=
\mid u_{max} - \max_{t \in(0,\infty)}\tilde{u}(t)\mid.
\label{eq:volterrapopulationerror}
\end{align}

\begin{figure}[!htbp]
    \centering
    
    \begin{subfigure}[b]{0.49\textwidth}
         \centering
         \includegraphics[width=\textwidth]{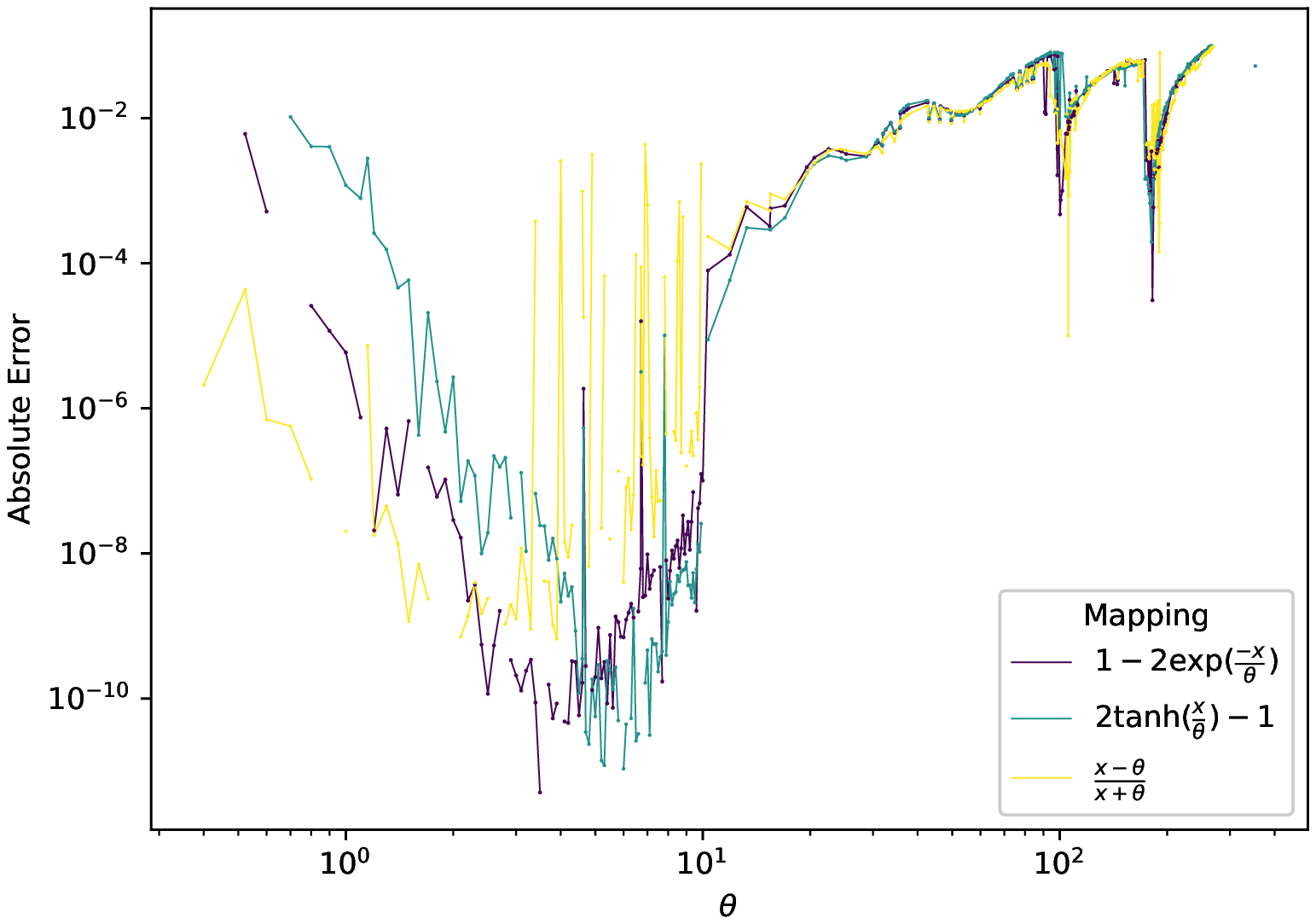}
         \caption{Legendre rational functions}
     \end{subfigure}
     \begin{subfigure}[b]{0.49\textwidth}
         \centering
         \includegraphics[width=\textwidth]{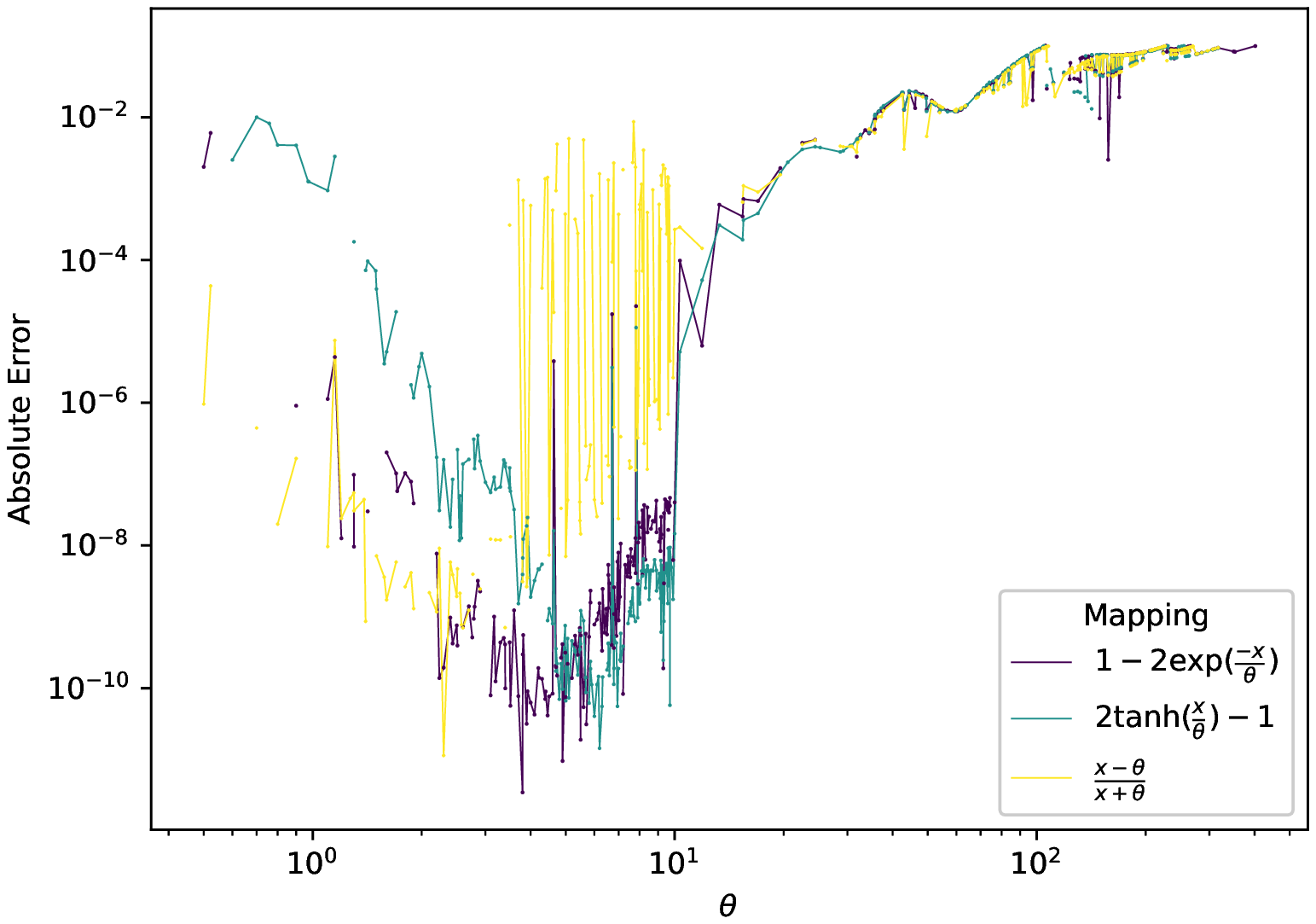}
         \caption{Chebyshev rational functions}
     \end{subfigure}
     \caption{The effect of length scale parameter $\theta$ on a large domain for Volterra's population model with $m=25$ and $\kappa=0.5$.}
     \label{fig:largevolterra}
\end{figure}

\begin{figure}[!htbp]
    \centering
    \includegraphics[width=0.7\textwidth]{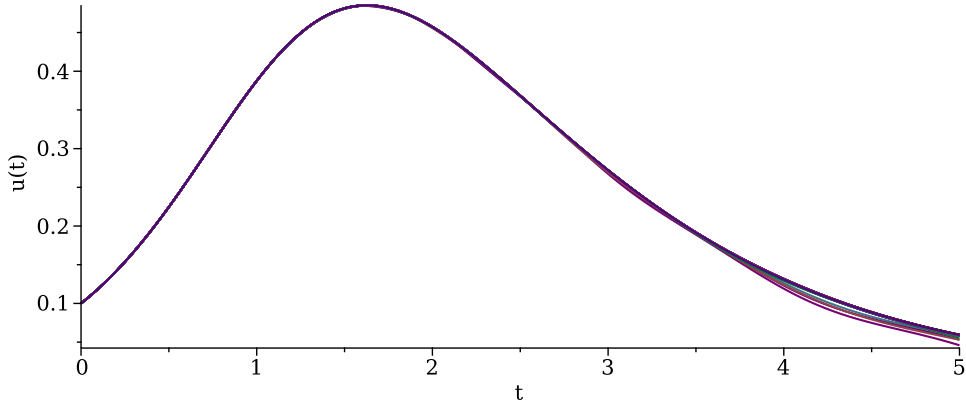}
    \caption{A set of learned solutions with different values for length scale $\theta$ for Volterra's population model with $\kappa=0.5$.}
    \label{fig:variousvolterra}
\end{figure}

To find a reasonable range for the length scale parameter and the effect of the other hyperparameters, we first ran a sensitivity analysis on these hyperparameters on a large domain. The results are reported in the figure \ref{fig:largevolterra}. It can be seen that the large values for the length scale do not yield a good approximation. The maximum reasonable value for this task is about $10$. In addition, the basis functions do not impose any significant difference on this parameter. Moreover, the nonlinear mappings can affect the computational procedure of the optimization problem \eqref{eq:opt}. This issue has resulted a discontinuity in these figures. Figure \ref{fig:variousvolterra} plots some of the successfully learned approximated solutions. It is seen that the equation may not be accurately simulated using improper hyperparameters. After the sensitivity analysis, we simulate this problem with five different most used values for the non-dimensional parameter $\kappa$ using the grid and random search algorithms. As reported in tables \ref{table:gridvalues} and \ref{table:randvalues}, the interval $(0,10]$ is chosen for the length scale parameter. Tables \ref{table:volterra-population-grid} and \ref{table:volterra-population-random} report the best-obtained hyperparameters. From there, it can be inferred that algebraic mapping is the best choice for small values of $\kappa$ while for the larger values, the exponential mapping could obtain better approximations. Moreover, the kernel function and its nonlinear mapping can result in a bit different accuracy.

To show the effectiveness of the proposed algorithm, we have compared the absolute error found by various authors in table \ref{table:volterra-population-comparison}. Rational Legendre (RL) and Rational Chebyshev (RC) pseudospectral method \cite{dehghan2015rational}, Sinc collocation method (SCM) \cite{parand2011collocation}, and Fractional Rational Legendre (FRL) \cite{parand2021parallel} are compared to each other. The length scale values used by Parand et al. \cite{parand2021parallel} are optimized in a similar manner proposed by Boyd \cite{boyd2001chebyshev}. 


\begin{table}[!htbp]
\centering
\resizebox{\textwidth}{!}{\begin{tabular}{@{}ccccccc@{}}
\toprule
$\kappa$ &    Kernel &                                       Mapping & $\theta$ &            Exact \cite{tebeest1997classroom} &      Approximate &     Error \\
\midrule
    $0.02$ &  Legendre &  $\displaystyle\nicefrac{(x-\theta)}{(x+\theta)}$ &  $0.10000$ & $0.92342717207022$ & $0.92342711545307$ & $5.662\times 10^{-08}$ \\
    $0.04$ &  Legendre &  $\displaystyle\nicefrac{(x-\theta)}{(x+\theta)}$ &  $0.70000$ & $0.87371998300000$ & $0.87371998000000$ & $3.090\times 10^{-09}$ \\
    $0.10$ &  Legendre & $\displaystyle1-2\exp(\nicefrac{-x}{\theta})$ &  $1.00000$ & $0.76974149100000$ & $0.76974149100000$ & $7.890\times 10^{-12}$ \\
    $0.20$ & Chebyshev & $\displaystyle1-2\exp(\nicefrac{-x}{\theta})$ &  $1.90000$ & $0.65905038200000$ & $0.65905038200000$ & $4.660\times 10^{-11}$ \\
    $0.50$ & Chebyshev & $\displaystyle1-2\exp(\nicefrac{-x}{\theta})$ &  $3.80000$ & $0.48519029140942$ & $0.48519029141289$ & $3.468\times 10^{-12}$ \\
\bottomrule
\end{tabular}}
\caption{The obtained results for Volterra population equation using a grid search algorithm.}
\label{table:volterra-population-grid}
\end{table}

\begin{table}[!htbp]
\centering
\resizebox{\textwidth}{!}{\begin{tabular}{ccccccc}
\toprule
$\kappa$ &    Kernel &                                        Mapping &          $\theta$ &               Exact \cite{tebeest1997classroom} &         Approximate &     Error \\
\midrule
    $0.02$ & Chebyshev &   $\displaystyle\nicefrac{(x-\theta)}{(x+\theta)}$ & $0.539501186666072$ &    $0.92342717207022$ &    $0.92342733722160$ & $1.652\times 10^{-07}$ \\
    $0.04$ & Chebyshev &   $\displaystyle\nicefrac{(x-\theta)}{(x+\theta)}$ & $0.318328463774207$ &    $0.87371998315400$ &    $0.87371998508417$ & $1.930\times 10^{-09}$ \\
    $0.10$ & Chebyshev &  $\displaystyle1-2\exp(\nicefrac{-x}{\theta})$ & $1.626117351946306$ &    $0.76974149070060$ &    $0.76974149073275$ & $3.216\times 10^{-11}$ \\
    $0.20$ & Chebyshev &  $\displaystyle1-2\exp(\nicefrac{-x}{\theta})$ & $2.510838579760311$ &    $0.65905038155232$ &    $0.65905038153414$ & $1.818\times 10^{-11}$ \\
    $0.50$ &  Legendre & $\displaystyle2\tanh(\nicefrac{x}{\theta})-1$ & $6.797026768536748$ &    $0.48519029140942$ &    $0.48519029141142$ & $2.007\times 10^{-12}$ \\
\bottomrule
\end{tabular}
}
\caption{The obtained results for Volterra population equation using a random search algorithm.}
\label{table:volterra-population-random}
\end{table}

\begin{table}[!htbp]
\centering
\begin{tabular}{@{}llllll@{}}
\toprule
$\kappa$ & RL \cite{dehghan2015rational} & RC \cite{dehghan2015rational} & FRL \cite{parand2021parallel} & SCM \cite{parand2011collocation} & Presented Method \\ 
$m$      &    $50$  &    $50$ &    $40$  &    $35$ &    $40$ \\
\midrule
$0.02$ & $3.72\times 10^{-07}$ & $7.51\times 10^{-07}$ & $6.33\times 10^{-06}$ & $7.00\times 10^{-08}$ & $5.66\times 10^{-08}$ \\
$0.04$ & $1.43\times 10^{-08}$ & $5.27\times 10^{-08}$ & $6.75\times 10^{-08}$ & $3.00\times 10^{-08}$ & $1.93\times 10^{-09}$ \\
$0.10$ & $1.07\times 10^{-10}$ & $2.13\times 10^{-10}$ & $2.73\times 10^{-08}$ & $8.00\times 10^{-08}$ & $7.89\times 10^{-12}$ \\
$0.20$ & $3.53\times 10^{-11}$ & $2.33\times 10^{-10}$ & $8.57\times 10^{-10}$ & $3.30\times 10^{-07}$ & $1.82\times 10^{-11}$ \\
$0.50$ & $2.44\times 10^{-09}$ & $4.87\times 10^{-09}$ & $2.73\times 10^{-10}$ & $3.40\times 10^{-07}$ & $2.01\times 10^{-12}$ \\ \bottomrule
\end{tabular}
\caption{A comparison among some mathematical methods solved Volterra's population model on the semi-infinite domain.}
\label{table:volterra-population-comparison}
\end{table}

\subsection{Kidder equation}
The unsteady isothermal flow of a gas through a micro-nano porous medium can be expressed as a nonlinear differential equation \cite{parand2018generalized}. This equation which is defined on a semi-infinite domain is modeled as:
\begin{equation}
\begin{aligned}
    &u''(x) + \frac{2x}{\sqrt{1- \kappa u(x) }} u'(x) = 0,\\
    &u(0) = 1, u(\infty)=0.
\label{eq:kidder}
\end{aligned}
\end{equation}
The initial slope of the approximated solution is an essential measure of the accuracy of the problem. Up to now, no exact solution or initial slope has been found for this problem. However, some researchers developed advanced techniques to find an accurate solution. In this research, we use the values obtained by Parand et al. \cite{parand2018accurate} as an exact approximation. This paper has obtained the exact initial slope up to $38$ digits of precision in a machine or software that supports arbitrary precision arithmetic. Here, we are just focusing on the problem of choosing the best hyperparameters. Thus, we compare our results on a small number of basis and an almost similar number of basis functions.
Furthermore, they have utilized the Quasi-Linearization Method (QLM) to convert the problem of approximating solutions for nonlinear Ordinary Differential Equations (ODEs) into a sequence of dependent linear ODEs. This increases the computational complexity to the number of QLM iterations. Here we have just solved the original nonlinear problem, which is more computationally efficient.

\begin{figure}[!htbp]
    \centering
    
    \begin{subfigure}[b]{0.49\textwidth}
         \centering
         \includegraphics[width=\textwidth]{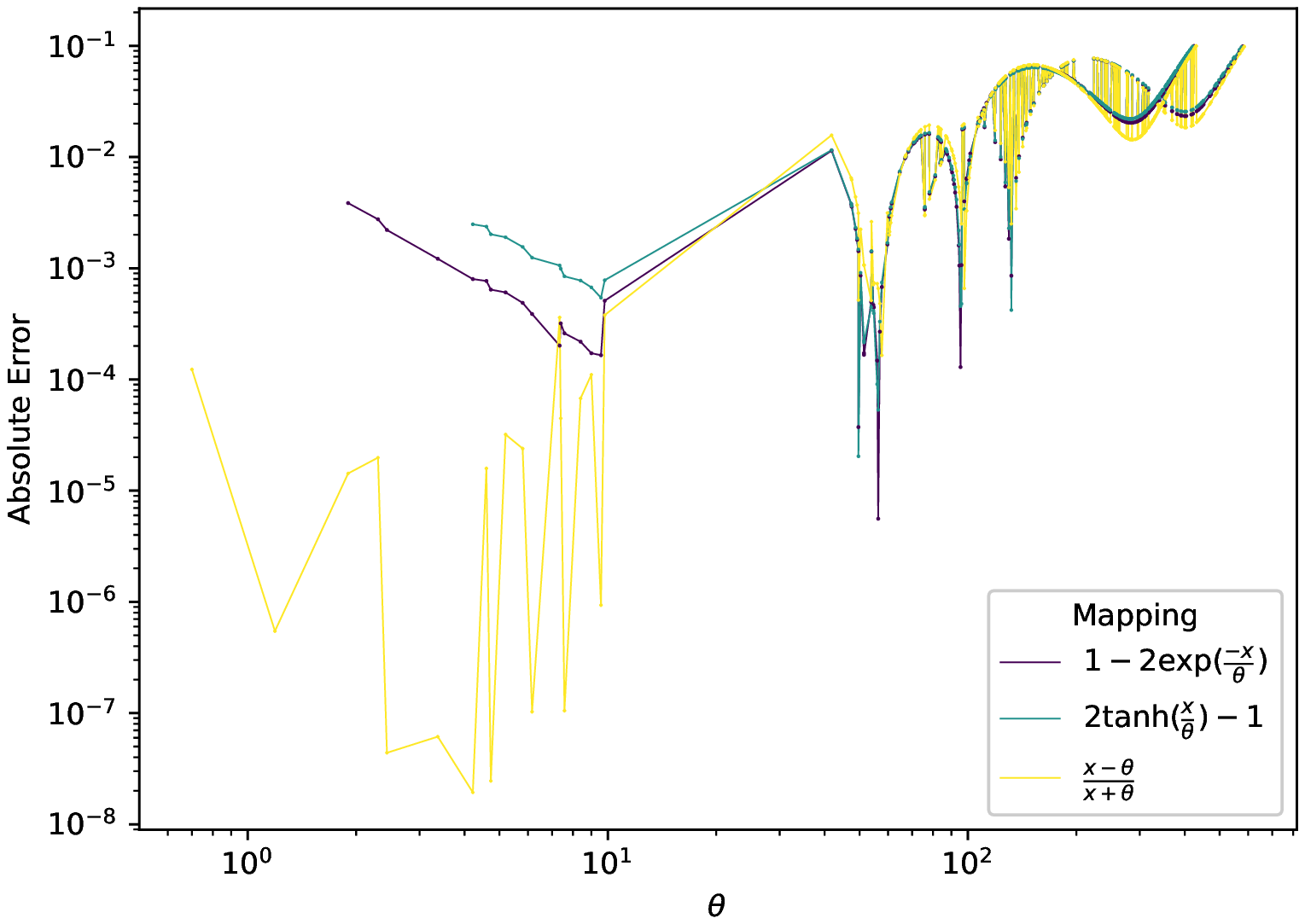}
         \caption{Legendre rational functions}
     \end{subfigure}
     \begin{subfigure}[b]{0.49\textwidth}
         \centering
         \includegraphics[width=\textwidth]{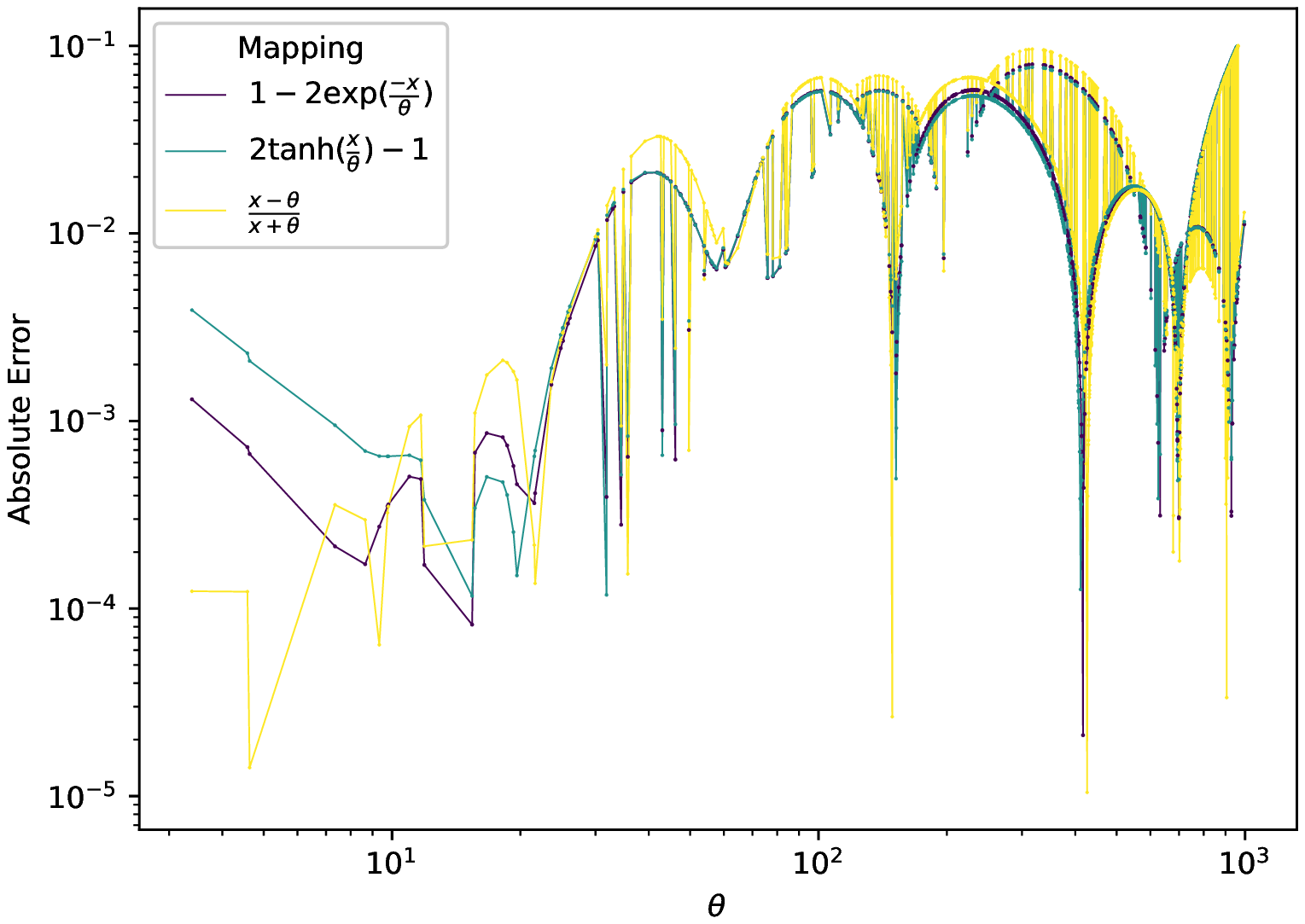}
         \caption{Chebyshev rational functions}
     \end{subfigure}
     \caption{The effect of length scale parameter $\theta$ on a large domain for Kidder equation with $m=25$}
     \label{fig:largekidder}
\end{figure}

Same as the previous example, we first analyze the effect of the hyperparameters. Figure \ref{fig:largekidder} demonstrates the absolute error obtained by different sets of hyperparameters. It is seen that the Legendre kernel can reach better results in comparison to the Chebyshev functions. Furthermore, some of the hyperparameters lead to numerical issues with the Legendre basis functions, and hence the plot is discontinuous. Likewise in the previous example, the interval $(0, 10]$ for the length scale contains the best results, therefore we focus on this range in the next experiments. Tables \ref{table:unsteady-gas-grid} and \ref{table:unsteady-gas-random} reported the best results obtained by the grid and random search, respectively. From there, it can be seen that the Legendre with algebraic mapping is the best-obtained hyperparameters in all of the configurations for different parameters. Furthermore, the random search algorithm overcame the grid search in all of the experiments.

\begin{figure}[!htbp]
    \centering
    \includegraphics[width=0.7\textwidth]{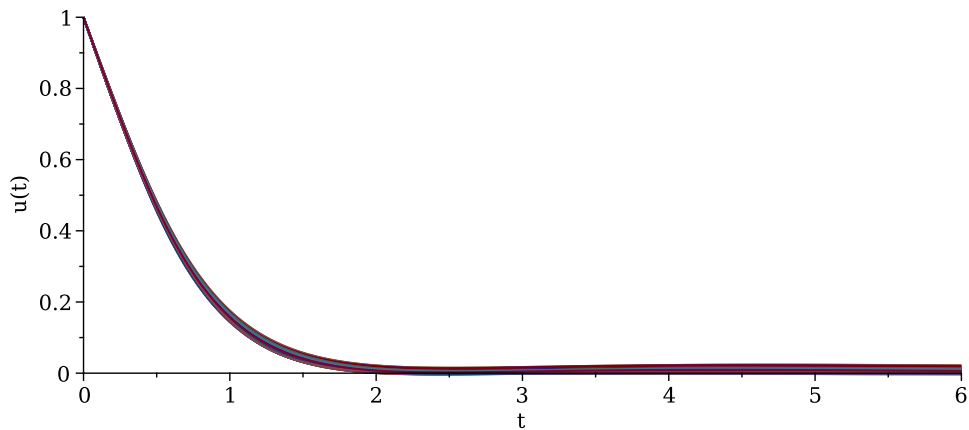}
    \caption{A set of learned solutions with different values for length scale $\theta$ for the Kidder equation with $\kappa = 0.5$.}
    \label{fig:variousunsteady}
\end{figure}

A comparison is made between the presented method and other related works in table \ref{table:kidder-comparison}. Some of the simulated approximations are plotted in the figure \ref{fig:variousunsteady}.
\begin{table}[!htbp]
\centering
\resizebox{\textwidth}{!}{\begin{tabular}{ccccccc}
\toprule
$\kappa$ &   Kernel &                                        Mapping & $\theta$ &             Exact \cite{parand2018accurate} &       Approximate &     Error \\
\midrule
    $0.10$ & Legendre &   $\displaystyle\nicefrac{(x-\theta)}{(x+\theta)}$ &  $4.80000$ & $-1.13900720617830$ & $-1.13900783881046$ & $6.326\times 10^{-07}$ \\
    $0.30$ & Legendre &   $\displaystyle\nicefrac{(x-\theta)}{(x+\theta)}$ &  $4.80000$ & $-1.16294145829591$ & $-1.16294198206589$ & $5.238\times 10^{-07}$ \\
    $0.50$ & Legendre &   $\displaystyle\nicefrac{(x-\theta)}{(x+\theta)}$ &  $4.80000$ & $-1.19179064971942$ & $-1.19179149343671$ & $8.437\times 10^{-07}$ \\
    % $0.70$ & Legendre & $\displaystyle-1+2\tanh(\nicefrac{x}{\theta})$ &  $4.00000$ & $-1.16294145829591$ & $-1.22562521115434$ & $6.268\times 10^{-02}$ \\
    $0.90$ & Legendre &   $\displaystyle\nicefrac{(x-\theta)}{(x+\theta)}$ &  $6.20000$ & $-1.28188132220336$ & $-1.28188069301075$ & $6.292\times 10^{-07}$ \\
\bottomrule
\end{tabular}}
\caption{The obtained results for the Kidder equation using a grid search algorithm.}
\label{table:unsteady-gas-grid}
\end{table}
\begin{table}[!htbp]
\centering
\resizebox{\textwidth}{!}{\begin{tabular}{ccccccc}
\toprule
$\kappa$ &   Kernel &                                      Mapping &          $\theta$ &             Exact \cite{parand2018accurate} &       Approximate &     Error \\
\midrule
    $0.10$ & Legendre & $\displaystyle\nicefrac{(x-\theta)}{(x+\theta)}$ & $0.975404049994095$ & $-1.13900720617830$ & $-1.13900677991113$ & $4.263\times 10^{-07}$ \\
    $0.30$ & Legendre & $\displaystyle\nicefrac{(x-\theta)}{(x+\theta)}$ & $1.576130816775483$ & $-1.16294145829591$ & $-1.16294118714719$ & $2.711\times 10^{-07}$ \\
    % $0.30$ & Legendre & $\displaystyle\nicefrac{(x-\theta)}{(x+\theta)}$ & $1.576130816775483$ & $-1.16294145829591$ & $-1.16294118713600$ & $2.712\times 10^{-07}$ \\
    $0.50$ & Legendre & $\displaystyle\nicefrac{(x-\theta)}{(x+\theta)}$ & $2.760250769985784$ & $-1.19179064971942$ & $-1.19179065703565$ & $7.316\times 10^{-09}$ \\
    $0.90$ & Legendre & $\displaystyle\nicefrac{(x-\theta)}{(x+\theta)}$ & $4.693906410582058$ & $-1.28188132220336$ & $-1.28188182747538$ & $5.053\times 10^{-07}$ \\
\bottomrule
\end{tabular}}
\caption{The obtained result for the Kidder equation using a random search algorithm.}
\label{table:unsteady-gas-random}
\end{table}
\begin{table}[!htbp]
\centering
\resizebox{\textwidth}{!}{\begin{tabular}{@{}cccccccc@{}}
\toprule
$\kappa$ & CPA \cite{iacono2015kidder} & Sinc \cite{rezaei2011numerical} & RL \cite{rezaei2011numerical} & HFC \cite{rad2011numerical} & Bessel \cite{parand2014solving} & RJ \cite{parand2017new} & Presented method \\ 
$m$ & $-$ & 32 & 32 & 20 & 26 & 20 & 25 \\ \midrule
$0.1$ & $2.94\times 10^{-04}$ & $-$ & $-$ & $-$ & $2.61\times 10^{-06}$ & $1.29\times 10^{-04}$ & $4.26\times 10^{-07}$\\
$0.5$ & $9.15\times 10^{-04}$ & $3.10\times 10^{-03}$ & $3.10\times 10^{-03}$ & $7.98\times 10^{-03}$ & $3.42\times 10^{-05}$ & $1.17\times 10^{-04}$ & $7.32\times 10^{-09}$\\
$0.9$ & $1.72\times 10^{-02}$ & $-$ & $-$ & $-$ & $3.65\times 10^{-06}$ & $1.02\times 10^{-04}$ & $5.05\times 10^{-07}$ \\ \bottomrule
\end{tabular}}
\caption{Comparison among various approximate initial slopes for the Kidder equation \eqref{eq:kidder}. The values reported in \cite{parand2018generalized} are assumed as the exact values.}
\label{table:kidder-comparison}
\end{table}

\section{Conclusion}
This paper developed two machine learning techniques for increasing the accuracy of the numerical simulation of functional equations. The presented algorithms \ref{alg:proposedgrid} and \ref{alg:proposedrand} are general tuning algorithms capable of hyperparameter optimization of various mathematical approaches such as spectral methods, RBFs, and wavelets. However, in this research, we have just focused on the spectral method for approximating the solution of nonlinear differential equations on the semi-infinite domain. To do so, we configured the search space to the hyperparameters of rational Jacobi functions, including basis functions, nonlinear rational mappings, and the length scale parameter. Finally, in the numerical results, various experiments have been conducted to measure the search capability of the proposed algorithms. We discussed the role of the length scale parameter on the stability and convergence of the method. In addition, some comparisons among related works were carried out to show the superiority of these algorithms over traditional mathematical procedures. This property, along with the small computation complexity handled with parallel programming approaches, made this process efficient and easy to use for other researchers. However, modern gradient-free global optimization techniques, such as Bayesian optimization and Tree-structured Parzen estimator, can be developed to get better approximations.

% \printbibliography

\bibliography{references.bib} 
\bibliographystyle{ieeetr}

\end{document}